\documentstyle[12pt]{article}
\newlength{\abstractwidth}
\renewcommand{\thefootnote}{\fnsymbol{footnote}}
\renewcommand{\thanks}[1]{\footnote{#1}} 
\newcommand{\starttext}{ \setcounter{footnote}{0}
\renewcommand{\thefootnote}{\arabic{footnote}}}

\newcommand{\be}{\begin{equation}}
\newcommand{\bea}{\begin{eqnarray}}
\newcommand{\eea}{\end{eqnarray}} \newcommand{\ee}{\end{equation}}
 
 \def\ba{\begin{eqnarray}}
\def\ea{\end{eqnarray}}

\def\tr{{\rm tr}}

\def\log{\,{\rm log}\,}
\def\exp{\,{\rm exp}\,}

\def\ge{\geq}
\def\le{\leq}

\def\p{\partial}

\def\[{{\bf [}}
\def\]{{\bf ]}}

\def\ddbar{i\p\bar\p}

\def\mathbb{\bf}

\def\eqref{\ref}

\begin{document}
\starttext \baselineskip=18pt \setcounter{footnote}{0}
\newtheorem{theorem}{Theorem}
\newtheorem{lemma}{Lemma}
\newtheorem{corollary}{Corollary}
\newtheorem{definition}{Definition}
\newtheorem{conjecture}{Conjecture}
\newtheorem{proposition}{Proposition}
\renewcommand{\eqref}[1]{(\ref{#1})}


\begin{center}
{\Large \bf On $L^\infty$ estimates for Monge-Amp\`ere and Hessian equations on nef classes
\footnote{Work supported in part by the National Science Foundation under grant DMS-1855947. }}

\medskip
\centerline{Bin Guo, Duong H. Phong, Freid Tong\footnote{F.T. is supported by Harvard's Center for Mathematical Sciences and Applications. }, and Chuwen Wang}

\medskip

\begin{abstract}

{\footnotesize The PDE approach developed earlier by the first three authors for $L^\infty$ estimates for fully non-linear equations on K\"ahler manifolds is shown to apply as well to Monge-Amp\`ere and Hessian equations on nef classes. 
In particular, one obtains a new proof of the estimates of Boucksom-Eyssidieux-Guedj-Zeriahi and Fu-Guo-Song for the Monge-Amp\`ere equation, together with their generalization to Hessian equations.}

\end{abstract}

\end{center}

\baselineskip=15pt
\setcounter{equation}{0}
\setcounter{footnote}{0}

\section{Introduction}
\setcounter{equation}{0}

The goal of this short note is to show that the PDE approach introduced in \cite{GPT, GPTa} for $L^\infty$ and Trudinger-type estimates for general classes fully non-linear equations on a compact K\"ahler manifold applies as well to Monge-Amp\`ere and Hessian equations on nef classes.

\smallskip
The key to the approach in \cite{GPT, GPTa} is an estimate of Trudinger-type, obtained by comparing the solution $\varphi$ of the given equation to the solution of an auxiliary Monge-Amp\`ere equation with the energy of the sublevel set function $-\varphi+s$ on the right hand side. We shall see that, in the present case of nef classes, the argument can still be made to work by replacing $\varphi$ by $\varphi-V$, where $V$ is the envelope of the nef class. Applied to the Monge-Amp\`ere equation, this gives a PDE proof of the estimates obtained earlier for nef classes by Boucksom-Eyssidieux-Guedj-Zeriahi \cite{BEGZ} and Fu-Guo-Song \cite{FGS}. The estimates which we obtain with this method applied to Hessian equations seem new.

\smallskip
We note that the use of an auxiliary Monge-Amp\`ere equation had been instrumental in the recent progress of Chen and Cheng \cite{CC} on the constant scalar curvature K\"ahler metrics problem. There the auxiliary equation involved the entropy, and not the energy of sublevel set functions as in our case. More generally, auxiliary equations have often been used in the theory of partial differential equations, notably by De Giorgi \cite{DeG}, and more recently by Dinew and Kolodziej \cite{DK, DDK} in their approach to H\"older estmates for the complex Monge-Amp\`ere equation.

\section{The Monge-Amp\`ere equation}
\setcounter{equation}{0}
\label{section 2}

We begin with the Monge-Amp\`ere equation.
Let $(X,\omega)$ be a compact K\"ahler manifold and $\chi$ be a closed $(1,1)$-form on $X$. We assume the cohomology class $[\chi]$ is nef and let $\nu\in \{0,1,\ldots, n\}$ be the numerical dimension of $[\chi]$, i.e.
$$\nu = \max\{k~|~ [\chi]^k\neq 0 \mbox{ in } H^{k,k}(X,\mathbb C)\}.$$
When $\nu = n$ we say the class $[\chi]$ is {\em big}.

Let $\hat \omega_t = \chi + t \omega$ for $t\in (0,1]$. The form $\hat\omega_t$ may not be positive but its class is K\"ahler. We consider the following family of complex Monge-Amp\`ere equations
\begin{equation}\label{eqn:MA}
(\hat \omega_t + \ddbar \varphi_t ) ^n = c_t e^{F} \omega^n, \quad\sup_X \varphi_t = 0
\end{equation}
where $c_t = [\hat \omega_t]^n = O(t^{n-\nu})$ is a normalizing constant and $F\in C^\infty(X)$ satisfies $\int_X e^F \omega^n = \int_X \omega^n$. This equation admits a unique smooth solution $\varphi_t$ by Yau's theorem \cite{Y}.

The form $\chi$ is not assumed to be semipositive, so the usual $L^\infty$ estimate of $\varphi_t$ may not hold \cite{K}. As in \cite{BEGZ, FGS}, we need to modify the solution $\varphi_t$ by an envelope $V_t$ of the class $[\hat \omega_t]$, defined as follows,
$$V_t = \sup\{v ~|~ v\in PSH(X,\hat \omega_t), \, v\le 0\}.$$
Then we have:

\begin{theorem}
\label{main}
Consider the equation (\ref{eqn:MA}), and assume that the cohomology class of $\chi$ is nef. For any $s>0$, let $\Omega_s = \{\varphi_t - V_t \le -s\}$ be the sub-level set of $\varphi_t - V_t$.

{\rm (a)} Then there are constants $C=C(n, \omega, \chi)>0$ and $\alpha_0 = \alpha_0(n,\omega,\chi)>0$ such that
\bea
\label{key}\int_{\Omega_s} \exp\big\{ \alpha_0 \Big( \frac{- (\varphi_t - V_t + s)}{A_s^{1/(1+n)}} \Big)^{\frac{n+1}{n}}    \big\}
\omega^n\le C\exp{(CE_t)},
\eea
where $A_s = \int_{\Omega_s} (-\varphi_t + V_t - s) e^F \omega^n$ and $E_t = \int_X (-\varphi_t + V_t)e^F \omega^n$.

{\rm (b)} Fix $p>n$.  There is a constant $C(n,p,\omega,\chi, \| e^F\|_{L^1(\log L)^p}   )$ so that for all $t\in (0,1]$, we have
\bea
0\le -\varphi_t+V_t\leq C(n,p,\omega,\chi, \| e^F\|_{L^1(\log L)^p}   ).
\eea
\end{theorem}

We remark that the estimates in Theorem \ref{main} continue to hold for a family of K\"ahler metrics (maybe with distinct complex structures) which satisfy a uniform $\alpha$-invariant type estimate. 

\medskip
\noindent
{\it Proof of Theorem \ref{main}}

We would like to find an auxiliary equation with smooth coefficients, so that its solvability can be guaranteed by Yau's theorem. For this, we 
need a lemma due to Berman \cite{B} on a smooth approximation for $V_t$ (see also Lemma \ref{lemma B1} below). Fix a time $t\in (0,1]$.

\begin{lemma}
\label{lemma Berman}
Let $u_\beta$ be the smooth solution to the complex Monge-Amp\`ere equation
$$(\hat \omega_t + \ddbar u_\beta)^n = e^{\beta u_\beta} \omega^n.$$ Then $u_\beta$ converges uniformly to $V_t$ as $\beta\to \infty$.
\end{lemma}

We remark that by \cite{CTW}, $V_t$ is a $C^{1,1}$ function on $X$, although this fact is not used in this note. We now return to the proof of Part (a) of Theorem \ref{main}.

\smallskip

We choose a sequence of positive functions $\tau_k: \mathbb R\to \mathbb R_+$ such that $\tau_k(x) $ decreases to $x\cdot\chi_{\mathbb R_+}(x)$ as $k\to\infty$. Fix a smooth function $u_\beta$ as in Lemma \ref{lemma Berman}. The $u_\beta$ depends on $t$, but for simplicity we omit the subscript $t$. We solve the following auxiliary Monge-Amp\`ere equation on $X$
\begin{equation}\label{eqn:auxMA}
(\hat \omega_t + \ddbar \psi_{t,k})^n = c_t \frac{\tau_k(-\varphi_t + u_\beta - s)}{A_{s,k,\beta}} e^F \omega^n,\quad \sup_X \psi_{t,k} = 0,
\end{equation}
where $A_{s,k,\beta} =  \int_X \tau_k(-\varphi_t + u_\beta - s) e^F \omega^n.$ Since $\psi_{t,k}\le V_t$ and $u_\beta$ converges uniformly to $V_t$. By taking $\beta$ large enough, we may assume $\psi_{t,k} < u_\beta+1$.

Define a function
$$\Phi = - \varepsilon ( -\psi_{t,k} + u_\beta + 1+ \Lambda)^{n/(n+1)} - (\varphi_t - u_\beta  + s),$$
with the constants
\begin{equation}
\varepsilon^{n+1} = A_{s,k,\beta}n^{-n} (n+1)^n,\quad \Lambda = n^{n+1} (n+1)^{-n-1} \varepsilon^{n+1}.
\end{equation}
As a smooth function on the compact manifold $X$, $\Phi$ must achieve its maximum at some $x_0\in X$. If $x_0\in X\backslash \Omega_s^\circ$, then
$$\Phi(x_0)\le - (\varphi_t - u_\beta  + s)\le - V_t + u_\beta \le \epsilon_\beta$$
where $\epsilon_\beta\to 0$ as $\beta\to \infty$. On the other hand, if $x_0\in \Omega_s^\circ$ we calculate ($\Delta_t$ denotes the Laplacian with respect to the metric $\omega_t = \hat \omega_t + \ddbar\varphi_t$)
\bea
0& \ge & \Delta_t \Phi(x_0) \nonumber\\
& = & -\varepsilon \frac{n}{n+1} (-\psi_{t,k} + u_\beta + \Lambda + 1) ^{-\frac{1}{n+1}} \tr_{\omega_t}( -\ddbar \psi_{t,k} + \ddbar u_\beta  )- \tr_{\omega_t} (\ddbar \varphi_t  - \ddbar u_\beta) \nonumber\\
 & & + \frac{n\varepsilon}{(n+1)^2} (-\psi_{t,k} + u_\beta + 1 + \Lambda)^{-\frac{n+2}{n+1}} \tr_{\omega_t} i \partial (\psi_{t,k} - u_\beta)\wedge \bar \partial (\psi_{t,k} - u_\beta)\nonumber\\
 & \ge  &  \frac{n \varepsilon}{n+1} (-\psi_{t,k} + u_\beta + \Lambda + 1) ^{-\frac{1}{n+1}} \tr_{\omega_t}( \hat\omega_{t,\psi_{t,k}} - \hat\omega_{t, u_\beta}  )-n + \tr_{\omega_t} \hat\omega_{t, u_\beta} \nonumber\\
 & \ge  & \frac{n \varepsilon}{n+1} (-\psi_{t,k} + u_\beta + \Lambda + 1) ^{-\frac{1}{n+1}} n \Big( \frac{\hat \omega_{t,\psi_{t,k}}^n}{\omega_t^n} \Big)^{1/n} - n + (1 - \frac{n \varepsilon}{n+1} (-\psi_{t,k} + u_\beta + \Lambda + 1) ^{-\frac{1}{n+1}} )\tr_{\omega_t} \hat \omega_{t,u_\beta}\nonumber\\
 & \ge &  \frac{n^2 \varepsilon}{n+1} (-\psi_{t,k} + u_\beta + \Lambda + 1) ^{-\frac{1}{n+1}}  
  \Big(\tau_k(-\varphi_t + u_\beta -s ) A_{s,k,\beta}^{-1} \Big)^{1/n}
  - n + (1 - \frac{n \varepsilon}{n+1} \Lambda ^{-\frac{1}{n+1}} )\tr_{\omega_t} \hat \omega_{t,u_\beta}\nonumber\\
   & \ge &  \frac{n^2 \varepsilon}{n+1} (-\psi_{t,k} + u_\beta + \Lambda + 1) ^{-\frac{1}{n+1}}
  (-\varphi_t + u_\beta -s )^{1/n} A_{s,k,\beta}^{-1/n}
  - n \nonumber.
\eea
Therefore, at $x_0\in \Omega_s^\circ$ 
$$-(\varphi_t - u_\beta +s ) \le (\frac{n\varepsilon}{n+1})^n A_{s,k,\beta} (-\psi_{t,k} + u_\beta + \Lambda + 1)^{\frac{n}{n+1}} = \varepsilon (-\psi_{t,k} + u_\beta + \Lambda + 1)^{\frac{n}{n+1}} , $$
i.e. $\Phi(x_0)\le 0$. Combining the two cases, we conclude that $\sup_X \Phi\le \epsilon_\beta \to 0$ as $\beta\to \infty.$ It then follows that on $\Omega_s$
$$( - \varphi_t + u_\beta - s) ^{\frac{n+1}{n}} \le  C_n A_{s, k,\beta} ^{1/n} (-\psi_{t,k} + u_\beta + 1 + A_{s,k,\beta}) + \epsilon_\beta^{(n+1)/n} $$
Letting $\beta\to \infty$ we have
$$( - \varphi_t + V_t - s) ^{\frac{n+1}{n}} \le  C_n A_{s, k} ^{1/n} (-\psi_{t,k} + V_t+ 1 + A_{s,k}), $$
where $A_{s,k} = \int_X \tau_k(-\varphi_t + V_t +s) e^F \omega^n$. Observe that by definition $V_t\le 0$ and by the $\alpha$-invariant estimate \cite{H, Ti}, there exists an $\alpha_0(n,\omega, \chi)$ such that
\begin{equation}\label{eqn:1}
\int_{\Omega_s} \exp\Big( \alpha_0  \frac{( - \varphi_t + V_t - s) ^{\frac{n+1}{n}}}{ A_{s,k}^{1/n}  }     \Big) \omega^n \le \int_{\Omega_s} \exp \Big( \alpha_0 C_n (-\psi_{t,k} + 1 + A_{s,k})   \Big) \omega^n  \le C e^{C A_{s,k}}.
\end{equation}
Letting $k\to \infty$, we obtain
$$\int_{\Omega_s} \exp\Big( \alpha_0  \frac{( - \varphi_t + V_t - s) ^{\frac{n+1}{n}}}{ A_{s}^{1/n}  }     \Big) \omega^n \ \le C e^{C A_{s}}.$$ Part (a) of Theorem \ref{main} is proved by noting that $A_s\le E_t$ for any $s>0$.

\medskip

Once Part (a) of Theorem \ref{main} has been proved, Part (b) can be proved by following closely the arguments in \cite{GPT}.

\smallskip

Fix $p>n$, and define $\eta: {\mathbb R}_+\to {\mathbb R}_+$ by $\eta(x) = (\log(1+x))^p$. Note that $\eta$ is a strictly increasing function with $\eta(0) = 0$, and let $\eta^{-1}$ be its inverse function. Denote
\bea
v: = \frac{\alpha_0}{2}\big( \frac{-\varphi_t + V_t- s}{A_s^{1/(n+1)}} \big)^{(n+1)/n}
\eea
then by
the generalized Young's inequality with respect to $\eta$, for any $z\in\Omega_s$, 
   \bea
  v(z)^p e^{F(z)} & \le &\nonumber \int_0^{\exp({F(z)})} \eta(x) dx + \int_0^{v(z)^p} \eta^{-1}(y) dy\\
    &\le&\nonumber \exp({F(z)}) (1+  |F(z)|)^p + C(p) \exp({2 v(z)})
    \eea
We integrate both sides in the inequality above over $z\in \Omega_s$, and get by Part (a), Theorem \ref{main} that 
   \bea
  \int_{\Omega_s} v(z)^p e^{F (z)} \omega^n & \le&\nonumber \int_{\Omega_s} e^{F} (1+  |F(z)|)^p \omega^n + \int_{\Omega_s} e^{2v(z)} \omega^n\\
  & \le&\nonumber \| e^{F} \|_{L^1(\log L)^p} + C + C e^{C {E}_t},
  \eea
  where the constant $C>0$ depends only on $n,\omega_X, \chi$. In view of the definition of $v$, this implies  
   \begin{equation}\label{eqn:final 1} \int_{\Omega_s} (-\varphi_t + V_t- s)^{\frac{(n+1)p}{n}} e^{F(z) } \omega^n \le 2^p\alpha_0^{-p} A_s^{\frac{p}{n}} \big(\| e^{F} \|_{L^1(\log L)^p} + C + C e^{C  {E}_t}\big).
   \end{equation}
From the definition of $A_s$, it follows from H\"older inequality that
\bea
     A_s  & = &\nonumber \int_{\Omega_s} (-\varphi_t + V_t - s) e^{F} \omega^n \\
     & \le&\nonumber  \Big(  \int_{\Omega_s} (-\varphi_t +V_t - s) ^{\frac{(n+1)p}{n}} e^{F}\omega^n\Big)^{\frac{n}{(n+1)p}} \cdot \Big( \int_{\Omega_s} e^{F} \omega^n \Big)^{1/q}\\
     &\le&\nonumber  A_s^{\frac {1}{n+1}}\Big(  2^p\alpha_0^{-p} \big(\| e^{F} \|_{L^1(\log L)^p} + C + C e^{C {E}_t}\big)\Big)^{\frac{n}{(n+1)p}} \cdot \Big( \int_{\Omega_s} e^{F} \omega^n \Big)^{1/q}
\eea
where $q>1$ satisfies  $\frac{n}{p(n+1)} + \frac{1}{q} = 1$, i.e. $q = \frac{p(n+1)}{p(n+1) - n}$. The inequality above yields 
\begin{equation}\label{eqn:1.11}
   A_s\le \Big( 2^p\beta_0^{-p} \big(\| e^{F} \|_{L^1(\log L)^p} + C + C e^{C{ E}_t}\big)\Big)^{1/p} \cdot \Big( \int_{\Omega_s} e^{F} \omega^n \Big)^{\frac{1+n}{qn}}.
\end{equation}
Observe that the exponent of the integral on the right hand of \eqref{eqn:1.11} satisfies
$$\frac{1+n}{qn} = \frac{pn + p - n}{pn} = 1+ \delta_0>1,$$
for $\delta_0: = \frac{p-n}{pn}>0$. For notation convenience, set 
\begin{equation}\label{eqn:b0}B_0 :=  \Big(  2^p\alpha_0^{-p} \big(\| e^{F} \|_{L^1(\log L)^p} + C + C e^{C  { E}_t}\big)\Big)^{1/p}.
\end{equation}
From (\eqref{eqn:1.11}) we then get 
\begin{equation}\label{eqn:1.14}
     A_s \le B_0 \Big( \int_{\Omega_s} e^{F} \omega^n \Big)^{1+\delta_0}.
   \end{equation}
If we define $\phi:{\mathbb R}_+ \to {\mathbb R}_+$ by $\phi(s) :=  \int_{\Omega_s} e^{F} \omega^n$ then (\eqref{eqn:1.14}) and the definition of $A_s$ imply that 
\begin{equation}\label{eqn:final 10}
r \phi(s+r) \le B_0 \phi(s)^{1+\delta_0},\quad \forall r\in [0,1] \mbox{ and } s\ge 0.
\end{equation}
$\phi$ is clearly nonincreasing and continuous, so a De Giorgi type iteration argument shows that there is some $S_\infty$ such that $\phi(s) = 0$ for any $s\ge S_\infty$. This finishes the proof of the $L^\infty$ estimate of $\varphi_t - V_t$, combining with a bound on $E_t$ by $\| e^F\|_{L^1(\log L)^1}$ which follows from Jensen's inequality (c.f. Lemma 6 in \cite{GPT}).  The proof of Theorem \ref{main} is complete.

\medskip
Finally, we note the recent advances in the theory of envelopes in \cite{GL1} and \cite{GL2}, which can provide an approach to $L^\infty$ estimates for Monge-Amp\`ere equations on Hermitian manifolds.

\section{Complex Hessian equations}
\setcounter{equation}{0}

We explain in this section how the proof of Theorem \ref{main} can be modified to give a similar result for a degenerate family of complex Hessian equations. With the same notations as above, we consider the $\sigma_k$-equations
\begin{equation}\label{eqn:sigma k}
(\hat \omega_t + \ddbar \varphi_t)^k \wedge \omega^{n-k} = c_t e^{F} \omega^n,\quad \sup_X \varphi_t = 0.
\end{equation}
Define the envelope corresponding to the $\Gamma_k$-cone 
$$\tilde V_{t, k} = \sup\{v | ~ v\in SH_k(X,\omega,\hat \omega_t)\cap C^2,\,v\le 0\}$$
where $v\in SH_k(X,\omega,\hat \omega_t)\cap C^2$ means that the eigenvalue vector of the linear transformation $\omega^{-1}\cdot (\hat \omega_t + \ddbar v)$ lies in the $\Gamma_k$-cone.

Let $$E_t(\varphi_t) = \int_X(-\varphi_t +\tilde V_{t, k}) e^{nF/k} \omega^n$$ be the entropy associated to the equation (\ref{eqn:sigma k}) as in \cite{GPT} and let $\bar E_t$ be an upper bound of $E_t(\varphi_t)$. Then the following $L^\infty$-estimate holds for the solution $\varphi_t$ to \eqref{eqn:sigma k}.
\begin{theorem}
Let $\varphi_t$ be the solution to \eqref{eqn:sigma k}, then there exists a constant depending on $\bar E_t$, $\| e^{\frac{n}{k} F} \|_{L^1(\log L)^p}$, $\frac{c_t}{[\hat\omega_t]^k [\omega]^{n-k}}$, $p>n$ such that
$$0\le -\varphi_t + \tilde V_{t, k}\le C.$$
\end{theorem}
This theorem can be derived using a similar argument as in Section \ref{section 2} with suitable modifications for $\sigma_k$ equations, c.f. \cite{GPT}, so we omit the details. The only novel ingredient is the smooth approximation of $\tilde V_{t, k}$, as in Lemma \ref{lemma Berman}. One can adapt the method in \cite{B} to derive this required approximation. For the convenience of the reader, we present a sketch of the proof. 
\begin{lemma}\label{lemma B1}
Fix $t\in (0,1]$. There exists a sequence of smooth functions $u_\beta \in SH_k(X,\omega,\hat \omega_t)$ converging uniformly to $\tilde V_{t,k}$ as $\beta\to \infty$.

\end{lemma}
\noindent{\em Proof}. Let $u_\beta \in SH_k(X,\omega,\hat \omega_t)$ be the solution to the $\sigma_k$-equations
\begin{equation}\label{eqn:sigma k 1}
(\hat \omega_t + \ddbar u_\beta)^k\wedge \omega^{n-k} = c_t e^{\beta u_\beta} \omega^n,
\end{equation}
which admits a unique smooth solution by \cite{DKa}. We claim that there is a constant $C_t>0$ such that
$$\sup_X | u_\beta - \tilde V_{t,k}  |\le \frac{C_t \log \beta}{\beta},$$ 
from which the lemma follows. 

By the maximum principle, at a maximum point of $u_\beta$, $\ddbar u_\beta \le 0$, so $\beta u_\beta \le \log \frac{\hat \omega_t^k\wedge\omega^{n-k}}{c_t \omega^n}\le C_t$, that is $u_\beta - \frac{C_t}{\beta}\le 0$. By the definition of $\tilde V_{t,k}$, it follows that \begin{equation}\label{eqn:test}u_\beta - \frac{C_t}{\beta} \le \tilde V_{t,k}. \end{equation}

On the other hand, we fix a smooth $u\le 0$ such that $\hat \omega_t + \ddbar u>0$. Such a $u$ exists because $[\hat \omega_t]$ is a K\"ahler class by assumption. For any $v\in SH_k(X,\omega,\hat \omega_t)\cap C^2$ with $v\le 0$, we consider the barrier function
$$\tilde u = \frac 1 {\beta} u + (1-\frac{1}{\beta}) v - \frac{C_t' \log \beta}{\beta}  $$
where $C_t'>0$ is a large constant to be determined. By direct calculation, we have
$$(\hat \omega_t + \ddbar \tilde u)^k\wedge \omega^{n-k}\ge \frac{1}{\beta^n} (\hat \omega_t + \ddbar u)^k\wedge \omega^{n-k} \ge e^{\beta \tilde u} \omega^n$$
where the last inequality holds if we choose $C_t'$ large enough so that
$$e^{-C_t'\log \beta} \le \frac{1}{\beta^k} \min_X  \frac{(\hat \omega_t + \ddbar u)^k\wedge \omega^{n-k} }{\omega^n}.$$ Therefore we get
$$(\hat \omega_t + \ddbar \tilde u)^k\wedge \omega^{n-k} \ge e^{\beta(\tilde u - u_\beta)} (\hat \omega_t + \ddbar u_\beta)^k\wedge \omega^{n-k}. $$
At the maximum point of $\tilde u - u_\beta$, $(\hat \omega_t + \ddbar \tilde u)^k\wedge \omega^{n-k} \le (\hat \omega_t + \ddbar u_\beta)^k\wedge \omega^{n-k}$. This shows that $\tilde u - u_\beta\le 0$ on $X$. Taking supremum over all such $v$'s in $\tilde u$, it follows that
$$(1-\frac 1 \beta) \tilde V_{t,k} \le u_\beta + \frac{C_t\log \beta}{\beta}.$$ The lemma follows from this and \eqref{eqn:test}.

\bigskip


\noindent Department of Mathematics \& Computer Science, Rutgers University, Newark, NJ 07102 

\noindent bguo@rutgers.edu,

\medskip

\noindent Department of Mathematics, Columbia University, New York, NY 10027 

\noindent phong@math.columbia.edu, wang.chuwen@columbia.edu

\medskip
\noindent Center for Mathematical Sciences and Applications, Harvard University, Cambridge, MA 02138

\noindent ftong@cmsa.fas.harvard.edu

\end{document}